\documentclass{article}
 
\usepackage{amssymb,amsmath,amsthm} 
\usepackage{mathrsfs} 
 
\newcommand{\cf}{\emph{cf.}~} 
 
\def\B{\mathscr B} 
\def\C{\mathscr C} 
\def\E{\mathcal E} 
\def\F{\mathbb F} 
\def\G{\mathcal G} 
\def\H{\mathcal H} 
\def\I{\mathscr I} 
\def\K{\mathscr K} 
\def\L{\mathcal L} 
 
\def\N{\mathbb N} 
\def\R{\mathbb R} 
\def\T{\mathscr T}

\def\Z{\mathbb Z} 
 
\def\<{\left\langle} 
\def\>{\right\rangle} 
\def\({\left(} 
\def\){\right)} 
\def\[{\left[} 
\def\]{\right]} 
\def\ltwo{\mathsf{L}^{\:\!\!2}} 
\def\e{\mathop{\mathrm{e}}\nolimits} 
\def\d{\mathrm{d}} 
\def\supp{\mathop{\mathrm{supp}}\nolimits}

\newtheorem{Theorem}{Theorem}[section] 
 
\newtheorem{Lemma}[Theorem]{Lemma} 
 
\newtheorem{Corollary}[Theorem]{Corollary} 
\newtheorem{Proposition}[Theorem]{Proposition} 
\newtheorem{Definition}[Theorem]{Definition}

\begin{document} 
  %--------------------------------------------------------------------------------------------------- 
  % Title 
  %--------------------------------------------------------------------------------------------------- 
  \title{\Large\textbf{Toeplitz algebras and spectral results for the one-dimensional Heisenberg model}}
   
  \author{Mondher Damak$^1$, Marius M\u antoiu$^2$\\ and Rafael Tiedra de Aldecoa$^3$} 
  \date{\small}
  \maketitle
\vspace{-1cm}

  \begin{quote}
    \emph{
	\begin{itemize}
      \item[$^1$] D\'epartement de math\'ematiques, Universit\'e de Sfax,
	    B.\,P. 802-3018, Sfax, Tunisia 
      \item[$^2$] Institute of Mathematics ``Simion Stoilow'' of the Romanian Academy,
	    P.\,O. Box 1-764, Bucharest, RO-014700 
      \item[$^3$] D\'epartement de math\'ematiques, Universit\'e de Paris XI, 91405 Orsay Cedex France 
      \item[] \emph{E-mails:} Mondher.Damak@fss.rnu.tn, Marius.Mantoiu@imar.ro,
	    rafael.tiedra@math.u-psud.fr 
    \end{itemize}
}
  \end{quote}

  %--------------------------------------------------------------------------------------------------- 
  % Abstract 
  %--------------------------------------------------------------------------------------------------- 
   
  \begin{abstract} 
	We determine the structure of the spectrum and obtain non-propagation estimates for a class of 
	Toeplitz operators acting on a subset of the lattice $\Z^N$. This class contains the Hamiltonian 
	of the one-dimensional Heisenberg model. 
  \end{abstract} 
   
  %--------------------------------------------------------------------------------------------------- 
  \section{Introduction} 
  \setcounter{equation}{0} 
  %--------------------------------------------------------------------------------------------------- 
  Spectral theory of Toeplitz operators and its connection with $C^*$-algebras is a vast topic. We
  only indicate \cite{Douglas72} as a textbook systematization of part of the early theory. 
  Since the appearence of Douglas' book, the theory has evolved by extension and abstraction in 
  many directions. The works closest to the present paper are the ones relating Toeplitz operators 
  and Toeplitz algebras to ordered groups \cite{Murphy87,Murphy91,Xu99,XC99}. 
 
  Our interest in the topic has been aroused by the remark \cite{Streater67,Streater74,Reed/Simon} 
  that the one-dimensional Heisenberg Hamiltonian $H$ of ferromagnetism can be written as a 
  direct sum $H=\oplus_{N\in\N}H_N$, where $H_N$ can be interpreted as the Laplace operator on 
  the subgraph 
  \begin{equation*} 
	\Z^N_<:=\{(x_1,\ldots,x_N)\in\Z^N\:\!:\:\!x_1<x_2<\dots<x_N\} 
  \end{equation*} 
  of the standard (Cayley) graph $\Z^N$. In Section \ref{section Heisenberg} we show that $H_N$ 
  (when suitably restricting to $\Z^N_<$) is equal to the sum of a Toeplitz operator and a 
  multiplicative potential, both belonging to the Toeplitz-like $C^*$-algebra $\T^<(\Z^N)$ 
   generated by the unilateral shifts on $\Z^N_<$. 
 
  Although $\T^<(\Z^N)$ is of none of the types thoroughly studied in the literature, its structure 
  is simple enough to suggest spectral results for the operators it contains (see Section 
  \ref{Toeplitz algebras}). Since the elements of $\T^<(\Z^N)$ can be written as direct integrals over 
  a torus, we are mainly concerned with spectral properties of the fibers (this presents a particular 
  interest in the case of the one-dimensional Heisenberg model, \cf \cite{Reed/Simon,Yafaev00}). We 
  express in Section \ref{The spectrum} the essential spectrum of the fibers as a union of spectra of 
  a family of subhamiltonians, improving part of the statements of \cite{Zoladek82}, which concern 
  a larger class of operators but imposes an unnecessary exponential decay condition. We 
  consider that our formalism and proofs are much more simple and natural than those of 
  \cite{Zoladek82}. 
  In Section \ref{Localization} we show the following type of result concerning the whole 
  Hamiltonian $H_N$. If $\kappa$ is a continuous real function with suitable support, then there exists 
  a natural family of multiplication operators $\{\chi_n\}_{n\in\N}$ for which 
  $\|\chi_n\kappa(H_N)\|$ is arbitrarily small if $n$ is large enough. This can be reformulated in 
  terms of the evolution group $\{\e^{-itH_N}\}_{t\in\R}$: at energies belonging to $\supp(\kappa)$, 
  the system gouverned by $H_N$ stays ``out of $\supp(\chi_n)$'' uniformly in time. Practically 
  $\supp(\kappa)$ must not intersect the spectrum of a certain subhamiltonian associated to some 
  ideal of $\T^<(\Z^N)$. For the Heisenberg model we put into evidence a nice interpretation involving 
  cluster properties of $N$-magnon states. 
 
  We think that the structure of $\T^<(\Z^N)$ may also be crucial for proving finer spectral and 
  scattering properties. The obtention of a Mourre estimate, which is a first step in this direction, 
  is a problem under review. 
   
  Actually the setting generalizes from $\Z^N$ to $\(\Z^m_{\rm lex}\)^N$, where $\Z^m_{\rm lex}$ 
  is the group $\Z^m$ ordered lexicographically. This relies on rather deep results of \cite{Murphy91} 
  on Toeplitz algebras associated to certain ordered groups. The statements of the present article can 
  be pushed to this more general case; they just require more involved notations. Unfortunately, this is 
  \emph{not} the right framework to study the $m$-dimensional Heisenberg model (another type of 
  Toeplitz algebra is needed), so we refrained from giving explicit detailed results in this situation; see
  however Remark 4.3.. 
 
  Let us finally mention that our treatment has few direct connections with previous work on the 
  spectral theory of Toeplitz operators. We were actually guided by the $C^*$-algebra approach to 
  spectral analysis for Schr\"odinger operators, as in
  \cite{ABG,Damak/Georgescu00,Georgescu/Iftimovici,Mantoiu02,AMP}. 
 
  %--------------------------------------------------------------------------------------------------- 
  \section{The one-dimensional Heisenberg model}\label{section Heisenberg} 
  \setcounter{equation}{0} 
  %--------------------------------------------------------------------------------------------------- 
  In order to justify the class of operators we study, we present here briefly and rather formally 
  the one-dimensional Heisenberg model. Further details may be found in 
  \cite{Streater67,Streater74,Reed/Simon}. 
 
  We consider the one-dimensional lattice $\Z$ with a spin-$\frac12$ 
attached at each vertex. Let 
  \begin{equation*} 
    \F(\Z):=\{\alpha:\Z\to\{0,1\}\:\!:\:\!\supp(\alpha)\ \text{is finite}\}\:\!, 
  \end{equation*} 
  and write $\{e^0,e^1\}:=\{(0,1),(1,0)\}$ for the canonical basis of the spin-$\frac12$ Hilbert 
  space $\mathbb C^2$. For any $\alpha\in\F(\Z)$ we denote by $e^\alpha$ the element 
  $\{e^{\alpha(x)}\}_{x\in \Z}$ of the direct product $\prod_{x\in \Z}\mathbb C^2_x$. We 
  distinguish the vector $e^{\alpha_0}$, where $\alpha_0(x):=0$ for all $x\in \Z$. Each element 
  $e^\alpha$ is interpreted as a pure state of the system of spins, and $e^{\alpha_0}$ as its 
  ground state with all spins pointing down. The Hilbert space $\L$ of the system (which is spanned 
  by the states with all but finitely many spins pointing down) is the incomplete tensor product 
  \cite[Sec. 2]{Streater67}, \cite[Sec. 2]{Streater74} 
  \begin{equation*} 
    \L:=\bigotimes_{x\in\Z}^{\alpha_0}\mathbb C^2_x\equiv\text{closed span} 
    \left\{e^\alpha\:\!:\:\!\alpha\in\F(\Z)\right\}. 
  \end{equation*} 
  The dynamics of the spins is given by the nearest-neighbour Heisenberg Hamiltonian 
  \begin{equation}\label{Heisenberg model} 
    L:=-\frac12\sum_{|x-y|=1}\Big\{a\big[\sigma_1^{(x)}\sigma_1^{(y)} 
    +\sigma_2^{(x)}\sigma_2^{(y)}\big]+b\big[\sigma_3^{(x)}\sigma_3^{(y)}-1\big]\Big\}\:\!. 
  \end{equation} 
  The operator $\sigma_j^{(x)}$ acts in $\L$ as the identity operator on each factor $\mathbb C^2_y$, 
  except on the component $\mathbb C^2_x$ where it acts as the Pauli matrix $\sigma_j$. The scalars 
  $a,b\in\R$ prescribe the anisotropy of the system. The case $a\ne b$ corresponds to the XXZ model, 
  whereas the standard Heisenberg model is obtained for $a=b$.  
 
  Let $\H:=\ell^2[\F(\Z)]$. Then the Hilbert spaces $\L$ and $\H$ are isomorphic due to the unitarity of 
  the mapping $\varsigma:\L\to\H$ sending $e^\alpha$ onto $\chi_{\{\alpha\}}$ for all $\alpha$ 
  ($\chi_{\{\alpha\}}$ stands for the characteristic function of the singleton $\alpha$). In particular, the 
  set $\F(\Z)$ may be considered as the configuration space for the system of spins. The Hamiltonian 
  $L$ is unitarily equivalent to a difference operator in $\H$. Given $\alpha\in\F(\Z)$ and $x\in\alpha$, 
  $y\notin\alpha$, we write $\alpha^y_x$ for the function of $\F(\Z)$ such that 
  $\supp(\alpha^y_x)=\supp(\alpha)\sqcup\{y\}\setminus\{x\}$. 
 
  \begin{Lemma}\label{difference operator} 
    For any $f\in\H$, $\alpha\in\F(\Z)$, one has the equality 
    \begin{equation*} 
      (\varsigma^{-1}L\varsigma f)(\alpha)=-2\sum_{|x-y|=1}\alpha(x)\[1-\alpha(y)\] 
      \big[af(\alpha^y_x)-bf(\alpha)\big]\:\!. 
    \end{equation*} 
  \end{Lemma} 
 
  \begin{proof} 
    The claim follows from a direct calculation using the properties of the Pauli matrices. 
  \end{proof} 
 
  For $N\in\N$, let $\F_N(\Z):=\{\alpha\in\F(\Z)\:\!:\:\!\supp(\alpha)\text{ has $N$ elements}\}$ and 
  set $\H_N:=\ell^2[\F_N(\Z)]$. Lemma \ref{difference operator} shows that the subspace $\H_N$ of $\H$ 
  is left invariant by $\varsigma^{-1}L\varsigma$. This is due to the fact that the Hamiltonian $L$ 
  commutes with the (magnon) number operator $\frac12\sum_{x\in X}(\sigma^{(x)}_3+1)$. Moreover it is 
  straightforward to show that the restriction $H_N:=(\varsigma^{-1}L\varsigma)\upharpoonright\H_N$ 
  is bounded and symmetric. 
 
  The operator $H_N$ has a more convenient form in another representation, which deserves the 
  introduction of some notations. Elements of $\Z^N$ are denoted generically by 
  $\xi\equiv(x_1,\dots,x_N)$, $\eta\equiv(y_1,\dots,y_N)$ or $\zeta\equiv(z_1,\dots,z_N)$, and $P^<$ 
  stands for the multiplication operator in $\ell^2(\Z^N)$ by the characteristic function 
  $\chi_{\Z^N_<}$. We shall often interpret this projection as an 
  operator from $\ell^2(\Z^N)$ to $\ell^2(\Z^N_<)$. For each $\eta\in\Z^N$ we define the unitary 
  ``bilateral shift'' $u_\eta$ in $\ell^2(\Z^N)$ by 
  \begin{equation*} 
  \(u_\eta f\)(\xi):=f(\xi-\eta)\:\!,\quad f\in\ell^2(\Z^N),~\xi\in\Z^N. 
  \end{equation*}  
  We also define the ``unilateral shift'' $v_\eta^<:=P^<u_\eta\upharpoonright\ell^2(\Z^N_<)$, which is a 
  partial isometry in $\ell^2(\Z^N_<)$. Thus, for any $\varphi\in \ell^1(\Z^N)$, one can consider the 
  convolution operator  
  \begin{equation*} 
  C_\varphi:=\sum_{\eta\in\Z^N}\varphi(\eta)u_\eta 
  \end{equation*} 
  in $\ell^2(\Z^N)$, and the Toeplitz operator  
  \begin{equation*} 
  T^<_\varphi:=\sum_{\eta\in\Z^N}\varphi(\eta)v_\eta^< 
  \end{equation*} 
  in $\ell^2(\Z^N_<)$. These operators are related by the formula 
  $T^<_\varphi=P^<C_\varphi\upharpoonright\ell^2(\Z^N_<)$. Let us also consider the projections 
  $q^<_\eta:=v^<_\eta(v^<_\eta)^*$ on the range of the partial isometries $v_\eta^<$. One sees easily 
  that $q^<_\eta$ is the multiplication operator by the characteristic function of the set 
  $\Z^N_<\cap(\Z^N_<+\eta)$. Finally, among the multiplication operators in $\ell^2(\Z^N_<)$, we 
  distinguish those of the form  
  \begin{equation*} 
  V^<_\varphi:=\sum_{\eta\in\Z^N}\varphi(\eta)q_\eta^<\;\!,\quad\varphi\in\ell^1(\Z^N)\:\!. 
  \end{equation*} 
  We identify now the Heisenberg Hamiltonian to an operator of the form above. The set $S$ stands for 
  the collection of vectors $\{s^\pm_i\}_{i=1}^N\subset\Z^N$ with components 
  $(s^\pm_i)_j:=\pm\delta_{ij}$. Notice that $S$ is a symmetric family of generators for the group $\Z^N$.  
 
  \begin{Proposition}\label{H is Toeplitz} 
    The Hamiltonian $H_N$ is unitarily equivalent to the operator $T_\varphi^<+V_\psi^<$, where 
    $\varphi:=-2a\chi_S$ and $\psi:=2b\chi_S$. 
  \end{Proposition} 
 
  \begin{proof}  
    Let $\phi:\Z^N_<\to\F_N(\Z)$ be the one-to-one map $\xi\mapsto\chi_{\{x_1,\dots,x_N\}}$. Let 
    $\Phi:\H_N\to\ell^2(\Z^N_<)$ be the unitary operator given by $\Phi(f):=f\circ\phi$ for any 
    $f\in\H_N$. Then one has for $g\in\ell^2(\Z^N_<)$, $\xi\in\Z^N_<$ 
    \begin{equation*} 
      (\Phi H\Phi^{-1}g)(\xi)=-2\sum_{j=1}^N\hspace{-3pt} 
      \sum_{\begin{array}{c}\vspace{-16pt}\\\scriptstyle y=x_j\pm1\vspace{-3pt} 
	  \\\scriptstyle y\notin\{x_1,\ldots,x_N\}\end{array}}\hspace{-12pt} 
      \left\{ag\[\phi^{-1}(\chi_{\{x_1,\ldots,x_{j-1},y,x_{j+1},\ldots,x_N\}})\]-bg(\xi)\right\}. 
    \end{equation*} 
    On another hand 
    \begin{equation*} 
      \phi^{-1}(\alpha)=(\min\alpha,\min(\alpha\setminus\{\min\alpha\}),\ldots) 
    \end{equation*} 
    for any $\alpha\in\F_N(\Z)$. Thus 
    \begin{equation*} 
      \Phi^{-1}(\chi_{\{x_1,\ldots,x_{j-1},y,x_{j+1},\ldots,x_N\}}) 
      =(x_1,\ldots,x_{j-1},y,x_{j+1},\ldots,x_N) 
    \end{equation*} 
    if $y=x_j\pm1$, $y\notin\{x_1,\ldots,x_N\}$ and $x_1<\ldots<x_N$. This implies that 
    \begin{align*} 
      (\Phi H\Phi^{-1}g)(\xi)&=-2\sum_{s\in S} 
      \chi_{\Z^N_<}(\xi-s)\[ag(\xi-s)-bg(\xi)\]\\ 
      &=-2\sum_{s\in S}[(av_s^<-bq_s^<)g](\xi)\\ 
      &=[(T_\varphi^<+V_\psi^<)g](\xi)\:\!.\qedhere 
    \end{align*} 
  \end{proof} 
 
  Proposition \ref{H is Toeplitz} is one of the motivations to study operators of the form 
  $T_\varphi^<+V_\psi^<$ for general functions $\varphi,\psi\in \ell^1(\Z^N)$. Actually, we can even 
  indicate a larger class of operators $T_\varphi^<+V_\psi^<$ which have an interpretation on their 
  own, outside Toeplitz theory. If $M\subset\Z^N$ a finite subset such that $\eta\in M$ implies
  $-\eta\in M$, then one can associate  to it a Cayley graph having $\Z^N$ as set of vertices by declaring
  that $\eta$ and $\zeta$ are connected iff $\eta-\zeta\in M$. The Laplace operator $\Delta^M$ of this
  Cayley graph is a convolution operator in $\ell^2(\Z^N)$, easy to understand when applying a Fourier
  transformation. On the other hand Laplacians $\Delta^M_E$ on subgraphs $E$ of this Cayley graph
  could be complicated objects. However it easy to show the identity
  $\Delta^M_{\Z^N_<}=T^<_{\chi_M}+V^<_{-\chi_M}$, which applies to the Heisenberg model in the
  case $M=S$. Thus the Laplacians $\Delta^M_{\Z^N_<}$ are subject of the spectral results that follow
  below. 
 
  %--------------------------------------------------------------------------------------------------- 
  \section{Toeplitz algebras}\label{Toeplitz algebras} 
  \setcounter{equation}{0} 
  %--------------------------------------------------------------------------------------------------- 
 
  In this section we collect some results on the Toeplitz-like algebra $\T^<(\Z^N)$. We first 
  introduce the appropriate abstract setting. 
 
  Let $X$ be a (discrete) abelian group and $E$ be a non-void subset of $X$. The projection 
  $P^E:\ell^2(X)\to\ell^2(E)$ is defined in $\ell^2(X)$ as the multiplication operator by the 
  characteristic function $\chi_E$. As before we introduce the unitary translation operators 
  $\{u_\eta\}_{\eta\in X}$ in $\ell^2(X)$ and the partial isometries 
  $\{v_\eta^E:=P^Eu_\eta\upharpoonright\ell ^2(E)\}_{\eta\in X}$ in $\ell^2(E)$. Once again, for 
  $\varphi\in\ell^1(X)$, the ``Toeplitz operator''  
  \begin{equation*} 
	T_\varphi^E:=\sum_{\eta\in X}\varphi(\eta)v_\eta^E 
  \end{equation*}  
  and the ``potential''  
  \begin{equation*} 
	V_\varphi^E:=\sum_{\eta\in X}\varphi(\eta)q_\eta^E 
	\equiv\sum_{\eta\in X}\varphi(\eta)v_\eta^E(v_\eta^E)^* 
  \end{equation*}  
  are available as operators in $\ell^2(E)$. 
 
  \begin{Definition}\label{Toeplitz algebra} 
    The $C^*$-algebra $\T^E(X)\subset\B[\ell^2(E)]$ generated by the family 
	$\{v_\eta^E\}_{\eta\in E}$ is called the \emph{Toeplitz algebra of the group $X$ with respect to 
	the subset $E$}. 
  \end{Definition} 
 
  Obviously $\T^E(X)$ contains all the operators of the form $T_\varphi^E+V_\psi^E$, 
  $\varphi,\psi\in\ell^1(\Z^N)$ (and many others). In fact $\T^E(X)$ is also generated by the family 
  $\{T^E_\varphi\:\!:\:\!\varphi\in\ell^1(X)\}$ as in the case of usual Toeplitz algebras. If $E=X$, 
  $\T^E(X)$ is equal to $C^*_{\rm r}(X)$, the reduced group $C^*$-algebra of $X$. Since $X$ is 
  abelian, we may identify it with $C^*(X)$, the envelopping $C^*$-algebra of the convolution Banach 
  $*$-algebra $\ell^1(X)$. The identification puts into correspondence $u_\eta$ with 
  $\delta_\eta:=\chi_{\{\eta\}}$. Due to the Fourier transform, $C^*(X)$ is isomorphic to the 
  $C^*$-algebra $C(\widehat X)$ of continuous complex functions on $\widehat X$, where 
  $\widehat X$ is the (compact, abelian) dual group of $X$. The conventions ``in dimension zero'' are 
  clear: The group is $X=\{0\}$, with subset $E=\{0\}$, all functions are scalars, and for 
  $\varphi\in\mathbb C$ one sets $T_\varphi=V_\varphi:=\varphi1$ in $\ell^2(X)\equiv\mathbb C$.  
 
  The algebras $\T^E(X)$ were mainly studied for $X$ an ordered group and $E$ its positive cone 
  \cite{Murphy87,Murphy91,Xu99,XC99}. The standard case is the ``classical'' Toeplitz algebra $\T^\N(\Z)\equiv\T$ 
  \cite{Douglas72} associated to the unilateral shift on $\ell^2(\N)$ (often presented in a Fourier 
  transformed realization). In our proofs the isomorphic algebra $\T^*:=\T^{\N^*}(\Z)$, 
  $\N^*:=\{1,2,\ldots\}$, will appear more naturally. The case $X=\mathbb Z^2$, $E=\mathbb N^2$ is studied in \cite{DH71}. The most relevant Toeplitz algebra for us is 
  $\T^<(\Z^N):=\T^{\Z^N_<}(\Z^N)$, which is \emph{not} of ordered type. We shall point out its 
  structure, which will be the main tool in analysing the operators $T^<_\varphi+V^<_\psi$. 
 
  The key facts are: 
  \begin{enumerate} 
	\item[(A)] If $\theta:X\to X'$ is a group isomorphism sending $E$ onto $E'$, then $\T^E(X)$ and 
	$\T^{E'}(X')$ are naturally isomorphic, the element $v_\eta^E$ of $\T^E(X)$ being sent onto 
	$v^{E'}_{\theta(\eta)}\in\T^{E'}(X')$.  
	\item[(B)] If $E_j$ is a subset of a group $X_j$, $j=1,\ldots,m$, then 
	$\T^{E_1\times\cdots\times E_m}(X)$ can be identified with the spatial tensor product 
	$C^*$-algebra $\bigotimes_{j=1}^m\T^{E_j}(X_j)$, $v^E_{(\eta_1,\ldots,\eta_m)}$ being identified 
	with $\otimes_{j=1}^mv^{E_j}_{\eta_j}$. 
  \end{enumerate} 
  Both isomorphisms are unitarily implemented, but this will not be used explicitly in the sequel. 
 
  Let $\theta:\Z^N\to\Z^N$ be the group automorphism defined by  
  \begin{equation*} 
	\theta(y_1,\ldots,y_N):=(y_1,y_2-y_1,\ldots,y_N-y_{N-1})\:\!, 
  \end{equation*}  
  with inverse  
  \begin{equation*} 
	\theta^{-1}(z_1,\ldots,z_N)=(z_1,z_1+z_2,\ldots,z_1+\dots+z_N)\:\!. 
  \end{equation*} 
  Obviously, $\E:=\theta(\Z^N_<)=\Z\times(\N^*)^{N-1}$. Thus, by (A), 
  $\T^<(\Z^N)$ and $\T^\E(\Z^N)$ are isomorphic, $v^<_\eta$ being sent onto $v^\E_{\theta(\eta)}$ for 
  any $\eta\in\Z^N$. Consequently this isomorphism sends $T^<_\varphi$ onto 
  $T^\E_{\varphi\circ\theta^{-1}}$ and $V^<_\psi$ onto $V^\E_{\psi\circ\theta^{-1}}$. In the next few 
  lines we only consider the case of $T^<_\varphi$. 
 
  By applying (B), one gets an isomorphism between $\T^\E(\Z^N)$ and 
  $\T^\Z(\Z)\otimes\[\T^{\N^*}(\Z)\]^{\otimes(N-1)}\equiv C^*(\Z)\otimes(\T^*)^{\otimes(N-1)}$, under which 
  $T^\E_{\varphi\circ\theta^{-1}}$ is transformed into  
  \begin{equation}\label{avatar} 
	\sum_{(z_1,\ldots,z_N)\in\Z^N}(\varphi\circ\theta^{-1})(z_1,\ldots,z_N)\, 
	\delta_{z_1}\otimes v^{\N^*}_{z_2}\otimes\cdots\otimes v^{\N^*}_{z_N}.
  \end{equation} 
  Now we apply the partial Fourier transform on the first variable. It maps $\delta_{z_1}\in C^*(\Z)$ onto
  $\e_{z_1}\in C(\mathbb T)$,  where $\mathbb T$ is equal to the interval $[0,1]$ (with $0$ identified to $1$) and 
  $\e_z(\tau):= \e^{-2\pi iz\tau}$ for all 
  $z\in\Z$, $\tau\in\mathbb T$. Namely we use $\mathscr F_1:=\mathscr F\otimes 1$, where 
  $(\mathscr Ff)(\tau):=\sum_{z\in\Z}\e_z(\tau)f(z)$  for all $\tau\in\mathbb T$, $f\in\ell^1(\Z)$. When 
  performing the sum over $z_1$, the operator (\ref{avatar}) becomes 
  \begin{equation}\label{avatar'} 
	\sum_{(z_2,\ldots,z_N)\in\Z^{N-1}}\[\mathscr F_1(\varphi\circ\theta^{-1})\] 
	(\:\!\cdot\:\!,z_2,\ldots,z_N)\,v^{\N^*}_{z_2}\otimes\cdots\otimes v^{\N^*}_{z_N}, 
  \end{equation} 
  which will be regarded as an element of 
  $C(\mathbb T)\otimes\T^{\otimes(N-1)}\equiv C\[\mathbb T;\T^{\otimes(N-1)}\]$. Similar arguments can be 
  carried on for $V_\psi^<$. When summing up we get the following lemma. For any $\rho\in\ell^1(\Z^N)$ and 
  $\tau\in\mathbb T$ we define $\mu(\tau)\rho\in\ell^1(\Z^{N-1})$ by 
  \begin{equation*} 
	\[\mu(\tau)\rho\](z_2,\ldots,z_N):=\[\mathscr F_1(\rho\circ\theta^{-1})\](\tau,z_2,\ldots,z_N)\:\!. 
  \end{equation*} 
 
  \begin{Lemma}\label{Toeplitz equivalences} 
	The $C^*$-algebras $\T^<(\Z^N)$ and $C(\mathbb T)\otimes\T^{\otimes(N-1)}$ are naturally 
	isomorphic. The isomorphism sends $T_\varphi+V_\psi$ onto the direct integral 
	\begin{equation*} 
		\int_\mathbb T^\oplus\d\tau\(T^{N-1}_{\mu(\tau)\varphi}+V^{N-1}_{\mu(0)\psi}\), 
	\end{equation*} 
	where the exponent of the Toeplitz operators refers to the subset $(\N^*)^{N-1}$ of the group 
	$\Z^{N-1}$. 
  \end{Lemma} 
  The presence of a direct integral is connected to the invariance of our operators under a natural action of $\Z$ 
  on $\Z^N_<$ by translations. For the Heisenberg model, this  can be 
  traced back to any of the earlier representations of the Hamiltonian. 
 
  As a rule, the spectral results will be stated only for operators of the form $T^<_\varphi+V^<_\psi$. By 
  introducing suitable notations, they could be extended to all the elements of $\T^{<}(\Z^N)$.  
 
  \begin{Corollary}\label{spectrum} 
    For any real functions $\varphi,\psi\in\ell ^1(\Z^N)$, one has 
    \begin{equation}\label{spectru} 
      \sigma_{\rm ess}\big(T^<_\varphi+V^<_\psi\big)=\sigma\big(T^<_\varphi+V^<_\psi\big) 
      =\bigcup_{\tau\in\mathbb T}\sigma\big(T^{N-1}_{\mu(\tau)\varphi}+V^{N-1}_{\mu(0)\psi}\big)\:\!. 
    \end{equation} 
  \end{Corollary} 
 
  \begin{proof} 
    The essential spectrum of an element $A$ of a $C^*$-algebra $\C$ composed of bounded operators in a 
    Hilbert space $\G$ coincides with the spectrum of its image in the quotient $\C/[\C\cap\K(\G)]$, 
    where $\K(\G)$ denotes the compact operators in $\G$. Thus the first equality follows from the obvious fact 
    that there is no non-trivial compact operator in $C(\mathbb T)\otimes\T^{\otimes(N-1)}\subset 
    \B[\ltwo(\mathbb T)]\otimes\B[\ell^2(\N^*)]^{\otimes(N-1)}\subset\B[\ell^2(\mathbb T\times(\N^*)^{N-1})]$. 
    For the second equality, we apply Lemma \ref{Toeplitz equivalences} and the discussion in \cite[Sec. 8.2.4]{ABG}
    on the spectrum of observables defined by a continuous family of selfadjoint operators. Note that the union on
    the r.h.s. is automatically closed.   
  \end{proof} 
 
  If $N=1$, then $T^<_\varphi$ is the convolution operator $C_\varphi$ and $V^<_\psi$ is 
  the multiplication operator by $\sum_{z\in\Z}\psi(z)=(\mathscr F\psi)(0)\in\R$. Thus the spectrum 
  (\ref{spectrum}) reduces to 
  \begin{equation*} 
    \sigma\big(T^<_\varphi+V^<_\psi\big)=(\mathscr F\varphi)(\mathbb T)+(\mathscr F\psi)(0) 
    =[\min(\mathscr F\varphi),\max(\mathscr F\varphi)]+(\mathscr F\psi)(0) 
  \end{equation*}  
  as it should be. 
 
  %--------------------------------------------------------------------------------------------------- 
  \section{The essential spectrum of the fiber Hamiltonians}\label{The spectrum} 
  \setcounter{equation}{0} 
  %--------------------------------------------------------------------------------------------------- 
 
  The fiber Hamiltonians $H(\tau):=T^{N-1}_{\mu(\tau)\varphi}+V^{N-1}_{\mu(0)\psi}$ can be 
  interpreted physically as ``energy operators at fixed quasi-momentum $\tau$''. We study now their 
  esssential spectrum. We use again the well known fact that, given a Hilbert space $\G$, the 
  essential spectrum of an operator $A\in\B(\G)$ coincides with the spectrum of its image in the 
  Calkin algebra $\B(\G)/\K(\G)$, where $\K(\G)$ is the ideal of compact operators. 
 
  The fibers $H(\tau)$ act in $\ell^2[(\N^*)^{N-1}]\simeq\ell^2(\N^*)^{\otimes(N-1)}$ and belong to 
  the $C^*$-algebra $(\T^*)^{\otimes(N-1)}$. So we are faced to the problem of understanding the 
  quotient of $(\T^*)^{\otimes(N-1)}$ by $\K[\ell^2(\N^*)^{\otimes(N-1)}]$, the latter being 
  identified to $(\K^*)^{\otimes(N-1)}$, where $\K^*:=\K[\ell^2(\N^*)]$. The discussion in \cite{DH71}
  is relevant here, especially in the case $N=2$. We start by proving a result in a more abstract setting. 
 
  \begin{Lemma}\label{quotients} 
	Let $\C_j$ be a nuclear $C^*$-subalgebra of $\B(\H_j)$, where $\H_j$ is a Hilbert space, $j=1,2$. 
	Let $\I_j$ be an ideal of $\C_j$ and let $\pi_j:\C_j\to\C_j/\I_j$ be the canonical $*$-morphism. 
	Then the mapping 
	\begin{align*} 
	  \overline\pi:\C_1\otimes\C_2&\to\[(\C_1/\I_1)\otimes\C_2\]\oplus\[\C_1\otimes(\C_2/\I_2)\]\\ 
	  A&\mapsto\{(\pi_1\otimes1)(A),(1\otimes\pi_2)(A)\} 
	\end{align*} 
	is a $*$-morphism, and $\ker(\overline\pi)=\I_1\otimes\I_2$. 
  \end{Lemma} 
 
  \begin{proof} 
	Since $\C_1$ and $\C_2$ are nuclear, the mappings $\pi_1\otimes1$ and $1\otimes\pi_2$ are 
	(surjective) $*$-morphisms with $\ker(\pi_1\otimes1)=\I_1\otimes\C_2$ and 
	$\ker(1\otimes\pi_2)=\C_1\otimes\I_2$ \cite[Thm. 6.5.2]{Murphy90}. From this it follows that 
	$\overline\pi$ is a $*$-morphism with 
    $\ker(\overline\pi)=(\I_1\otimes\C_2)\cap(\C_1\otimes\I_2)$. Since ideals in nuclear 
    $C^*$-algebras are nuclear, the triple $(\I_1,\C_2,\I_2)$ verifies the (right) slice map 
    conjecture \cite[Prop. 10]{Wassermann78}. Thus one gets the equality $\ker(\overline\pi)=\I_1\otimes\I_2$
    by an easy adaptation of the proof of \cite[Cor. 5]{Wassermann76}. 
  \end{proof}  
 
  We need two more notations. For any $j\in\{2,\ldots,N\}$, $\varrho\in\ell^1(\Z^{N-1})$ and 
  $\tau\in\mathbb T$, we define $\nu_j(\tau)\varrho\in\ell^1(\Z^{N-2})$ by 
  \begin{equation*} 
    [\nu_j(\tau)\varrho](z_2,\ldots,z_{j-1},z_{j+1},\ldots,z_N) 
    :=(\mathscr F_j\varrho)(z_2,\ldots,z_{j-1},\tau,z_{j+1},\ldots,z_N)\:\!, 
  \end{equation*} 
  $\mathscr F_j$ being the Fourier transformation in the $j^{\rm th}$ variable. Furthermore we 
  denote by $\Sigma_j(\tau,\tau')$ the spectrum of the Toeplitz operator (relative to the pair 
  $\{\Z^{N-2},(\N^*)^{N-2}\}$) 
  \begin{equation*} 
    T^{N-2}_{\nu_j(\tau')\mu(\tau)\varphi}+V^{N-2}_{\nu_j(0)\mu(0)\psi} 
  \end{equation*} 
  acting in $\ell^2[(\N^*)^{N-2}]$. 
 
  \begin{Theorem}\label{Toeplitz essential spectrum} 
    Let $\varphi,\psi\in\ell^1(\Z^N)$ be real functions and $\tau\in\mathbb T$. Then one has 
    \begin{equation}\label{decomposition on Z_nu (1)} 
      \sigma_{\rm ess}\big(T^{N-1}_{\mu(\tau)\varphi}+V^{N-1}_{\mu(0)\psi}\big) 
       =\bigcup_{j=2}^{N}\bigcup_{\tau'\in\mathbb T}\Sigma_j(\tau,\tau')\:\!. 
    \end{equation} 
  \end{Theorem} 
 
  \begin{proof} 
    By analogy to (the Fourier transformed version of) the isomorphism\newline
    ${\T/\K[\ell^2(\N)]\simeq C(\mathbb T)}$ \cite[Thm. 7.23]{Douglas72}, one has a canonical
    isomorphism\newline$\T^*/\K^*\simeq C(\mathbb T)$, uniquely defined by 
    the fact that for any $z\in\Z$ the operator $v_z^{\N^*}\in\T^*$ is sent onto the function 
    $\e_z\in C(\mathbb T)$. Therefore the class of operators $T_f^{\N^*}\in\T^*$, $f\in\ell^1(\Z)$, 
    must be sent onto the class of functions $\mathscr Ff\in C(\mathbb T)$. 
 
    The Toeplitz algebra $\T^*$ is nuclear, since it is the extension of the abelian quotient 
    $\T^*/\K^*\simeq C(\mathbb T)$ by the nuclear ideal $\K^*$. So we may use the analog of 
    Lemma \ref{quotients} for the $N-1$ factors $\C_2=\ldots=\C_N=\T^*$ and 
    $\I_2=\ldots=\I_N=\K^*$. We get an injective $*$-morphism 
    \begin{align*} 
      (\T^*)^{\otimes(N-1)}/(\K^*)^{\otimes(N-1)}\hookrightarrow 
      &\oplus_{j=2}^N\(\T_*^{\otimes(j-2)}\otimes C(\mathbb T)\otimes\T_*^{\otimes(N-j)}\)\\ 
      &\simeq\oplus_{j=2}^N\(C(\mathbb T)\otimes\T_*^{\otimes(N-2)}\), 
    \end{align*} 
    sending $T_{\mu(\tau)\varphi}+V_{\mu(0)\psi}$ onto the collection  
    \begin{equation*} 
      \Big\{\int_\mathbb T^\oplus\d\tau'\(T^{N-2}_{\nu_j(\tau')\mu(\tau)\varphi} 
      +V^{N-2}_{\nu_j(0)\mu(0)\psi}\)\Big\}_{j=2,\ldots,N}\:\!. 
    \end{equation*} 
    One concludes by using \cite[Sec. 8.2.4]{ABG}. 
  \end{proof} 
 
  As an example, if $N=2$, one has 
  \begin{equation*} 
    \sigma_{\rm ess}\big(T^{1}_{\mu(\tau)\varphi}+V^{1}_{\mu(0)\psi}\big) 
    =\left\{(\mathscr F_1\mathscr F_2\varphi)(\tau-\tau',\tau') 
    +(\mathscr F_1\mathscr F_2\psi)(0,0)\:\!:\:\!\tau'\in\mathbb T\right\}, 
  \end{equation*} 
  In the case of the Heisenberg model (see Proposition \ref{H is Toeplitz}), one has to take 
  $\varphi=-2a\chi_S$ and $\psi=2b\chi_S$. Since 
  $(\mathscr F_1\mathscr F_2\chi_S)(\tau-\tau',\tau')=2\cos(2\pi\tau')+2\cos[2\pi(\tau-\tau')]$, 
  the essential spectrum of the fiber operator coincides with the interval 
  $\{8b-4a\cos(2\pi\tau')-4a\cos[2\pi(\tau-\tau')]\:\!:\:\!\tau'\in[0,1]\}$.  
 
  {\bf Remark 4.3.} Results of this section (and also of the next one) can easily be generalized to
  Toeplitz-like operators $T^<_\varphi+V^<_\psi$ acting on the subset
  \begin{equation*} 
    E=\(\Z_{\rm lex}^m\)_<^N 
	:=\left\{(x_1,\ldots,x_N)\in(\Z_{\rm lex}^m)_<^N\:\!:\:\!x_1<x_2<\cdots<x_N\right\} 
  \end{equation*} 
  of $X=\(\Z_{\rm lex}^m\)^N$, where $\Z_{\rm lex}^m$ is the group $\Z^m$ with the lexicographic
  order. In analysing the structure of $\mathcal T^E(X)$ for this case, the Toeplitz algebra $\T$ should
  be replaced by $\T_m:=\T^{(\Z_{\rm lex}^m)_+}(\Z_{\rm lex}^m)$, where 
  $(\Z_{\rm lex}^m)_+:=\{x\in\Z_{\rm lex}^m\:\!:\:\!0\leq x\}$. The only important change concerns 
  the calculation of the quotients $\T_m/\K_m$, where $\K_m:=\K\{\ell^2[(\Z_{\rm lex}^m)_+]\}$. 
  Namely one has to call for the rather deep result \cite[Thm. 2.3]{Murphy91}, which implies that 
  $\T_m/\K_m$ and $\T_{m-1}\otimes C(\mathbb T)$ are $*$-isomorphic. 
 
  %--------------------------------------------------------------------------------------------------- 
  \section{Localization}\label{Localization} 
  \setcounter{equation}{0} 
  %--------------------------------------------------------------------------------------------------- 
 
  In the sequel we determine localization properties of the operators $T^<_\varphi+V^<_\psi$ by
  adapting to the Toeplitz algebra $\T^<(\Z^N)$ a technique developed in \cite{AMP} (see also
  \cite{Mantoiu03}) for crossed product $C^*$-algebras, with applications to Schr\"odinger operators
  in $\mathbb R^n$. 
 
  Let $H$ be a selfadjoint operator in $\ell^2(\Z^N_<)$ (or in some other $\ltwo$-space) and $\chi$ a non-trivial multiplication 
  operator (for example the characteristic function of a set having a strictly positive measure). If $\kappa$ is a continuous function 
  with support intersecting the spectrum of $H$, the operator $\chi\kappa(H)$ has no reason to be small in general. The unique a 
  priori bound would be  
  \begin{equation}\label{small} 
    \|\chi\kappa(H)\|\leq\|\chi\|_\infty\sup_{\lambda\in\sigma(H)}|\kappa(\lambda)|\:\!. 
  \end{equation} 
  We are going to correlate $\chi$ to $\kappa$ in such a way to make the norm small without asking any of the two factors on the 
  r.h.s. of (\ref{small}) to be small. 
 
  \begin{Theorem}\label{non-propagation} 
    Fix $j\in\{2,\ldots,N\}$ and let $\varphi,\psi$ be real elements of $\ell^1(\Z^N)$. For any $n\in\N$ set  
    \begin{equation*} 
      \Omega_j(n):=\left\{(y_1,\ldots,y_N)\in\Z^N_<\:\!:\:\!y_j-y_{j-1}\ge n\right\}. 
    \end{equation*} 
    Let $\kappa:\mathbb R\to\mathbb R$ be a continuous function with  
    \begin{equation*} 
      \supp\kappa\cap\big[\textstyle{\bigcup_{\tau,\tau'\in\mathbb T}}\Sigma_j(\tau,\tau')\big]=\varnothing. 
    \end{equation*} 
    Then for each $\varepsilon>0$ there exists $n_\varepsilon\in\N$ such that 
    \begin{equation}\label{nii} 
      \big\|\chi_{\Omega_j(n)}\kappa(T^<_\varphi+V^<_\psi)\big\|\leq\varepsilon 
    \end{equation} 
    for each $n\ge n_\varepsilon$. 
  \end{Theorem} 
 
  \begin{proof} 
    Denote by $T^\E_{\varphi\circ\theta^{-1}}+V^\E_{\psi\circ\theta^{-1}}$ the image of $T^<_\varphi+V^<_\psi$ through 
    the isomorphism $\T^<(\Z^N)\simeq C(\mathbb T)\otimes(\T^*)^{\otimes(N-1)}$, defined by the change of variables 
    $\theta$. One verifies easily that the estimate (\ref{nii}) is equivalent to 
    \begin{equation}\label{iin} 
      \big\|\chi_{\widetilde{\Omega_j}(n)} 
      \kappa\big(T^\E_{\varphi\circ\theta^{-1}}+V^\E_{\psi\circ\theta^{-1}}\big)\big\|\leq\varepsilon\:\!, 
    \end{equation} 
    where 
    $\widetilde{\Omega_j}(n):=\theta\[\Omega_j(n)\]=\left\{(z_1,\ldots,z_N)\in\Z\times(\N^*)^{(N-1)}\:\!:\:\!z_j\ge n\right\}$. 
    Moreover the operator $\kappa\big(T^\E_{\varphi\circ\theta^{-1}}+V^\E_{\psi\circ\theta^{-1}}\big)$ belongs to the ideal 
    \begin{equation*} 
      \I_j:=C(\mathbb T)\otimes(\T^*)^{\otimes(j-2)}\otimes\K^*\otimes(\T^*)^{\otimes(N-j)} 
    \end{equation*} 
    of $\mathscr C:=C(\mathbb T)\otimes(\T^*)^{\otimes(N-1)}$. Indeed the image of 
    $T^\E_{\varphi\circ\theta^{-1}}+V^\E_{\psi\circ\theta^{-1}}$ in the quotient 
    \begin{equation*} 
      \mathscr C/\I_j\simeq C(\mathbb T)\otimes(\T^*)^{\otimes(j-2)}\otimes C(\mathbb T)\otimes(\T^*)^{\otimes(N-j)} 
      \simeq C(\mathbb T)^{\otimes 2}\otimes(\T^*)^{\otimes(N-2)} 
    \end{equation*} 
    is  
    \begin{equation*} 
      \int_{\mathbb T^2}^\oplus\d\tau\:\!\d\tau' 
      \(T^{N-2}_{\nu_j(\tau')\mu(\tau)\varphi}+V^{N-2}_{\nu_j(0)\mu(0)\psi}\), 
    \end{equation*} 
    with spectrum $\bigcup_{\tau,\tau'\in\mathbb T}\Sigma_j(\tau,\tau') 
    \equiv\sigma_{\I_j}\big(T^\E_{\varphi\circ\theta^{-1}}+V^\E_{\psi\circ\theta^{-1}}\big)$. Thus, since 
    $\supp(\kappa)\cap\big[\bigcup_{\tau,\tau'\in\mathbb T}\Sigma_j(\tau,\tau')\big]=\varnothing$, it 
    follows by \cite[Lemma 1]{AMP} that 
    $\kappa\big(T^\E_{\varphi\circ\theta^{-1}}+V^\E_{\psi\circ\theta^{-1}}\big)\in\I_j$. 
    Now 
    \begin{equation*} 
      \chi_{\widetilde{\Omega_j}(n)}=1_{j-1}\otimes\chi_{\{z_j\ge n\}}\otimes1_{N-j}, 
    \end{equation*} 
    where $\chi_{\{z_j\ge n\}}$ converges strongly to $0$ in $\B[\ell^2(\N^*)]$ as $n\to\infty$.  Thus, by examining the 
    structure of $\mathcal J_j$, one gets 
    \begin{equation*} 
      \big\|\chi_{\widetilde{\Omega_j}(n)}\kappa\big(T^\E_{\varphi\circ\theta^{-1}}+V^\E_{\psi\circ\theta^{-1}}\big)\big\| 
      \longrightarrow 0 
    \end{equation*} 
    as $n\to\infty$.  
  \end{proof} 
 
  Let $H$ be a selfadjoint operator in a Hilbert space $\G$ with spectral measure $E^H$, and let $f\in\G$ be an 
  arbitrary vector. We call \emph{spectral support of $f$ with respect to $H$}, and write $\supp(f;H)$, for the smallest closed 
  subset $F$ of $\mathbb R$ such that $E^H(F)f=f$. Alternatively one can characterize $\supp(f;H)$ as follows: 
  \begin{equation*} 
    \lambda\notin\supp(f;H)~~\text{iff}~~\exists\epsilon>0~\text{such that} 
    ~E^H(\lambda-\epsilon,\lambda+\epsilon)f=0. 
  \end{equation*} 
  Obviously one has $\text{supp}(f;H)\subset\sigma(H)$. If $H$ is the Hamilton operator describing some quantum 
  system in $\G$, we say that $f$ is \emph{\it a state with energy in $\supp(f;H)$}. 
   
  \begin{Corollary}\label{Toeplitz non-propagation} 
    Let $j,\varphi,\psi$ and $\kappa$ be as in Theorem \ref{non-propagation}. Then, for any $\varepsilon>0$, there 
    exists $n_\varepsilon\in\N$ such that 
    \begin{equation*} 
      \big\|\chi_{\Omega_j(n)}\e^{-it\(T^<_\varphi+V^<_\psi\)}f\big\|\leq\varepsilon\|f\| 
    \end{equation*} 
    for all $n\ge n_\varepsilon$, $t\in\R$ and all $f\in\ell^2(\Z^N_<)$ satisfying 
    \begin{equation*} 
      \supp\big(f;T^<_\varphi+V^<_\psi\big) 
      \cap\big[\textstyle{\bigcup_{\tau,\tau'\in\mathbb T}}\Sigma_j(\tau,\tau')\big]=\varnothing. 
    \end{equation*} 
  \end{Corollary} 
 
  Corollary \ref{Toeplitz non-propagation} follows trivially from Theorem \ref{non-propagation}. We put it into evidence 
  for its physical interpretation in the case of the one-dimensional Heisenberg model: Intuitively, if $f$ is a normalized 
  initial state with energy outside $\cup_{\tau,\tau'}\Sigma_j(\tau,\tau')$, the decomposition of the system into two 
  clusters of spins pointing up, one ``at the left'', composed of $j-1$ elements, and the other one ``at the right'', composed 
  of $N-j+1$ elements, is highly unprobable \emph{uniformly in time} if the distance $n$ between the clusters is large 
  enough. 
 
  It is obvious that several variants are available. One can consider ideals smaller than $\I_j$, by collapsing more than 
  one factor to the ideal of compact operators in $\ell^2(\N^*)$. In this way, one gets more detailed clustering information 
  for the one-dimensional Heisenberg model, but for a priori smaller sets of energy values. The fiber Hamiltonians 
  $T^{N-1}_{\mu(\tau)\varphi}+V^{N-1}_{\mu(0)\psi}$ can be studied identically by using ideals in the $C^*$-algebra 
  $(\T^*)^{\otimes(N-1)}$. 
 
  From Corollary \ref{spectrum} and Theorem \ref{Toeplitz essential spectrum} we know that 
  \begin{equation*} 
    \bigcup_{j=2}^N\bigcup_{\tau,\tau'\in\mathbb T}\Sigma_j(\tau,\tau')=\bigcup_{\tau\in\mathbb T}\sigma_{\rm ess} 
    \Big(T^{N-1}_{\mu(\tau)\varphi}+V^{N-1}_{\mu(0)\psi}\Big)\subset\sigma\big(T^<_\varphi+V^<_\psi\big). 
  \end{equation*} 
  It does not seem easy to determine under which conditions there is room in the spectrum of 
  $T^<_\varphi+V^<_\psi$ outside $\bigcup_{\tau,\tau'\in\mathbb T}\Sigma_j(\tau,\tau')$ for a given $j$. However one 
  may expect it is often the case. 
   
  %--------------------------------------------------------------------------------------------------- 
  \section*{Acknowledgements} 
  %--------------------------------------------------------------------------------------------------- 
 
  We are grateful to R. Purice for his helpful remarks. M. M{\u a}ntoiu acknowledges partial support 
  from the contract CERES, 4-187/2004. R. Tiedra de Aldecoa is supported by the Swiss National Science Foundation.
  Part of this work has been completed while M. M{\u a}ntoiu visited the University of Geneva; he expresses 
  his gratitude to W. Amrein for his kind hospitality. 
 
  %--------------------------------------------------------------------------------------------------- 
 
  \providecommand{\bysame}{\leavevmode\hbox to3em{\hrulefill}\thinspace} 
  \providecommand{\MR}{\relax\ifhmode\unskip\space\fi MR } 
  \providecommand{\MRhref}[2]{% 
    \href{http://www.ams.org/mathscinet-getitem?mr=#1}{#2} 
  } 
  \providecommand{\href}[2]{#2} 
   
  %--------------------------------------------------------------------------------------------------- 
 
\end{document}